\documentstyle{article}

\newcommand{\tu}{{\tilde u}}
\newcommand{\tv}{{\tilde v}} 
\newcommand{\hu}{{\hat u}}
\newcommand{\hv}{{\hat v}}
\newcommand{\n}{{{\cal N}_-}}
\newcommand{\g}{{\cal G}} 
\newcommand{\qed}{{$\Box $}}

\newcommand{\HH}{{\widetilde H}}

 \newcommand{\dfrac}[2]{{\frac{\textstyle #1}{\textstyle #2} }}

	 \newcommand{\Ri}{Riemannian }

	\newcommand{\sR}{subRiemannian }
	\newcommand{\R }{I \!\! R}

	\newtheorem{proposition}{Proposition}
	\newtheorem{theorem}{Theorem}
	\newtheorem{lemma}{Lemma}

	\title{\sc  Chaotic
	Geodesics in Carnot Groups}

	\author{
	Richard Montgomery \\ 
        email:  rmont@cats.ucsc.edu \\
        Mathematics Dept. UCSC, Santa Cruz, CA 95064, USA,\\
        Mikhail Shapiro, \\
                e-mail: mshapiro@math.kth.se \\
Dept. of Mathematics, KTH, 10044, Stockholm, Sweden, \\
                Alexander Stolin, \\
                 e-mail: astolin@math.chalmers.se \\
Dept. of Mathematics, University of G\"oteborg, \\
41296, G\"oteborg, Sweden.}

\begin{document}

	\maketitle

	{\sc abstract}

	Graded nilpotent Lie groups, or {\sc Carnot Groups}
	are to subRiemannian geometry as Euclidean
	spaces are to Riemannian geometry.  They are
	the metric tangent cones for this geometry.
	Hoping that the analogy between subRiemannian
	and Riemannian geometry is a strong
	one, one might conjecture that the  subRiemannian
	geodesic flow on any Carnot group is 
	completely integrable.    
	We  prove this conjecture is false by showing that
	that the
	subRiemannian geodesic flow is not 
	algebraically completely integrable
	in the case of the  group whose Lie algebra $\n$ consists
	of 4 by 4 nilpotent triangular matrices. 
We use this to prove that the  
  centralizer for the corresponding 
  quadratic ``quantum'' Hamiltonian $H$
in the universal enveloping algebra of $\n$
is ``as small as possible''.    
	\section{Introduction}

	Geometry would be in a sorry state
if the Euclidean geodesic flows were not completely
	integrable -- in other words, if we did  not have  
	explicit algebraic descriptions of 
	 straight lines in Euclidean space.  Riemannian
	geometry, being infinitesimally  Euclidean,  
makes frequent use of these explicit descriptions.
For example the Euclidean lines play  a role in the exponential map
and in Jacobi fields.   

	SubRiemannian geometries, also called Carnot-Caratheodory
	geometries, are not infinitesimally Euclidean.
	Rather they are, at typical points, infinitesimally
	modelled by Carnot groups.  We will review these
	geometries and groups, and the relation between them
	momentarily.  {\bf The point of this note is to show that
	the Carnot geodesic flows need not be integrable.
	We do this by giving an example of a `chaotic' Carnot
	geodesic flow.}  

	A subRiemannian geometry consists of  a nonintegrable
	subbundle (distribution) $V$ of the tangent bundle $T$ of
	a manifold, together with a fiber-inner product on this
	bundle.  These geometries arise, among other places,
	as the limits of Riemannian geometries.  In such a limit
	we penalize  curves for moving transverse to $V$ 
	so that in the limit any curve
	not tangent to $V$ has infinite length. 
	The distance between two points in
	a \sR manifold is defined as in \Ri geometry:
	it is the infimum of the lengths of all
	absolutely continuous paths connecting the two points.
	A theorem attributed to Chow asserts that
	if $V$ generates  $T$ under repeated Lie
	brackets, and if the manifold is connected
	then this distance function is everywhere finite.
Equivalently, and two points
can be connected by a curve tangent
to $V$.   
	In this manner, every \sR manifold becomes a metric
	space.

	By a Carnot group we mean a  simply connected Lie
	group  $G$ whose Lie algebra $\g$ is finite-dimensional,
	nilpotent, and graded with the degree 1 part
	generating the algebra and endowed with an
	inner product.  
	Specifically:  
	$$\g = V_1 \oplus V_2 \oplus \ldots V_r$$
	as a vector space.  The Lie bracket satisfies
	$$[V_i, V_j] \subset V_{i+j}$$
	where  $V_s = 0$ for $s > r$.  Also:  
	$$V_{i+1} = [V_1, V_i]$$
	and 
	$V_1$ is an inner-product space.
	We may think of $V= V_1$ as a left-invariant
	distribution on the group $G$.  Its
	inner product then gives $G$
	a subRiemannian geometry.

	Given a distribution 
	$V \subset T$, we can, at typical points
	$q$ of $Q$,  obtain
	a graded nilpotent Lie algebra.
	In order to do this, let 
	 $V$  also stand for the 
	the sheaf of  smooth  vector fields whose values lie in
	$V$.   Form  
	$$V ^{2} =  [ V , V] $$
	$$V ^{3} =  [V , V ^{2}]$$ 
	$$\vdots$$
	where the brackets denote Lie brackets of vector
	fields.
	(Exercise:  Show that, as sheaves:  $V^j \subset
	V^{j+1}$.) 
	We  assume that   $V$ is bracket
	generating,  which means that in a neighborhood of any
	point there is an integer $r$ such that  $V ^r =  T$.
	(This is the hypothesis of Chow's theorem.)
	 A
	point $q$ of the manifold is called {\sc regular} if the 
	sheaves $V^r$ correspond to vector bundles near $q$. This
	means that if we  evaluate each of these spaces of vector
	fields at a point $q \in Q$ thus obtaining a flag of
	subspaces:  \begin{equation} V_q \subset V^{2}_q \subset
	V^{3}_q \subset \ldots \subset V ^r (q) = T_q Q
	\label{eq:flag}
	\end{equation}
	then the integers $dim(V^j (q))$ are constant
	in some neighborhood of $q$. 
	For a generic distribution, generic points are
	regular.  
	Associate to the filtration $V  \subset V ^2
	\subset \ldots \subset T$  of sheaves its corresponding
	graded object $$ Gr(V , T)_q = V \oplus V_2 \oplus \ldots
	V_i \ldots \oplus V_r $$ where
	$V_j = V ^j  / V ^{j-1} $
	is the quotient sheaf.
	 Because  the Lie bracket
	of vector fields $X, Y$ satisfies
	$[X, fY] = f[X,Y] \hbox{ mod } X, Y$
	where $f$ is a function, it   induces 
	 bilinear maps
	$ : V_j \otimes V_k \to V_{j+k}$.  Putting
	these maps together defines, at any regular point, 
	 a Lie algebra structure on
	$Gr_q = Gr(V, T)(q) = V(q) \oplus V_2 (q) \oplus \ldots
	V_r (q)$. The subspace
	$V_q = V_1(q)$ of $Gr(V, T)_q$ is the original 
	k-plane field at that point, and
	Lie-generates $Gr_q$.   It follows
	that $\g = Gr(V, T)_q$ is
	the Lie algebra of a Carnot group.
	If the distribution $V$ comes
	with an inner product, then this generating
	subspace inherits it.  Consequently
	$G$ comes with a canonical left-invariant
	\sR structure.  This $G$
	is called the {\sc nilpotentization}
	of the \sR structure at the regular point $q$.

Gromov, using the idea of the Hausdorff
	limits of a family of metric
	spaces,  showed how to define a 
	``tangent space'' to any point of any metric space.
(See Gromov et al  \cite{Gr2}.) 
This limiting space,
called the metric tangent cone,
 only exists for `nice' metric
spaces.    
	For a  \Ri manifold it is 
	the usual tangent space with its Euclidean
structure.  The metric
	tangent cone also exists for  \sR metrics. A theorem of  
	of Mitchell \cite{Mitchell} asserts that
	 at a regular point it is the  nilpotentization.
	We urge the reader to consult Gromov \cite{Gr2}, \cite{Gr1}
	 and Bellaiche \cite{Bell}.
 
	The nilpotentization is the closest object we
	have in  \sR geometry to the Euclidean tangent
	space in \Ri geometry.  
	The match is not perfect but it is the best thing we
	have.

	\subsection{The geodesic flow}
	 
	The  data of a \sR geometry can be re-encoded as  a
	fiber-quadratic non-negative form $H: T^* Q \to \R$
	on the cotangent bundle $T^* Q$.  
	The kernel $\{ H = 0 \}$ of $H$ is the
	annihilator of the distribution $V$. Upon
	polarization $H$ becomes a bilinear non-negative form, and
	thus a symmetric map $g:T^* Q \to TQ$.  The image of
	this map is $V$ and the  map satisfive
	$\langle g(p), v \rangle _q =  p(v) $ for any $p \in T^Q*
	_q, v \in V_q, q \in Q$, where $\langle \cdot, \cdot
	\rangle$ is the inner product on $V$.     
	 The  Hamiltonian flow associated to $H$ generates curves
	in the cotangent bundle whose projections to Q
	are  subRiemannian geodesics.  By a \sR
geodesic we mean a curve in Q with the property
that the length of any sufficiently
short subarc of the curve equals the
\sR distance between the endpoints of
this arc.  Such curves are necessarily 
	necessarily tangent to $V$.

	To write down $H$ explicitly, pick
	any local orthonormal frame $\{X_1, \ldots , X_k\}$
	for the distribution $V$.
	Now the $X_i$ can be viewed as
	fiber-linear functions on the
	cotangent bundle $T^*$ and hence
	their squares $X_i^2$ are fiber-quadratic functions.
	The Hamiltonian is:
	$$H =  { 1 \over 2} (
	X_1 ^2 + X_2 ^2 +... X_k ^2 ).$$

	{\sc remark.}  Unlike \Ri geometry, there may be \sR
	geodesics which are not the projections
	of these solutions to Hamilton's equations in $T^* Q$. 
	See \cite{me1}, \cite{me2}, {me3}, \cite{LiuSuss}.
	But `most' geodesics are obtained as 
	projections of these solutions.

	In view of the analogies between \Ri
	and \sR geometries, the question naturally arises:
	{\bf Is geodesic flow on a Carnot group always
	integrable?} The answer is
	``yes'' for two step nilpotent groups.
	(The flow on the Heisenberg group has
	been integrated in many places.)  
	 The purpose of this
	note is to provide an example of a 3 step Carnot group 
	whose \sR geodesic flow is {\bf not} algebraically
	completely integrable.  Roughly speaking, this means that
	there is no uniform algebraic description of
	its ``straight lines''.
	 We expect  nonintegrability to
	hold for generic r step
	nilpotent graded groups , $r > 2$. 

In order to proceed we need to describe
how the geodesic flow for a Lie group
can be pushed down to a hamiltonian flow
on the dual of its Lie algebra.  We will also
need to recall the definition of ``complete
integrability''.

	In the particular case where the \sR geometry
	is that of a Carnot group $G$, then the 
	frame $X_i$ for $V = V_1$ can be realized by left-invariant
	vector fields.  We may identify
	the space of left-invariant vector fields
	with the Lie algebra. Thus
	$X_i \in \g$.  Then 
	$H$ becomes identified with a fiber-quadratic
	function on the dual $\g ^*$
	of the Lie algebra of $G$. 
	We recall that the dual of any Lie-algebra
	has a Poisson structure, the so-called ``Lie-Poisson
	structure''.  This can be defined by insisting that
	$$\{X_i, X_j \} = -[X_i, X_j]$$
	where we identify elements $X_i$ of the Lie algebra
	$\g$ with linear functions on its dual
	$\g^*$.  Thus
	$H$ induces a Hamiltonian  vector-field on
	$\g^*$.      

	Geometrically, what we are doing by studying
	this Hamiltonian flow on $\g^*$ is studying the  `Poisson
	reduction' of the \sR geodesic flow on $T^*G$. 
	The function $H$ on $T^*G$ is left-invariant, and
	hence so is its \sR geodesic flow.  The
	 vector field defining this
	flow can then  be pushed down to the quotient space
	$(T^*G)/G$ of the cotangent bundle by the
	left $G$ action, thus defining
	the ``Poisson-reduced'' flow. 
	Now $(T^*G)/G = \g^*$ in a natural way
	and when we  push down the Hamitonian
	vector field for $H$ we obtain the  one discussed
	in the  previous paragraph on $\g^*$.

	We recall the definition
	of completely integrable.
	A Hamiltonian $H$ (or its flow)
	on a {\bf symplectic} manifold
	of dimension $2n$ 
	is called {\sc completely integrable}
	if we can find $n$ functions
	$f_1, \ldots, f_n$
	which are almost everywhere 
	functionally independent
	($df_1 \wedge \ldots df_n \ne 0$),
	 which
	Poisson-commute ($\{f_i, f_j\} = 0$), and such that 
	$H$ can be expressed as a function of
	them ($H = h(f_1, \ldots , f_n)$).   If the flows
	of the $f_i$ are complete then their
	common level sets $\{f_1 = c_1, \ldots f_n = c_n\}$
	are, for typical constants $c_i$,   diffeomorphic
	to the  quotient of
	$R^n$ by a lattice, and on each
	such level, the flow of $H$ is linear
	up in the covering space $\R^n$.  
	The diffeomorphism is provided by action
	angle coordinate.  

	We are interested in whether the \sR geodesic flow on $T^*G$, 
$G$ a Carnot group, is
 complete integrable. 
	In order to proceed    
	{\sc we will  assume
	that if the flow `upstairs' on
	$T^*G$ is completely integrable,
	then so is the flow `downstairs'
	on $\g^*$.}  The converse 
 is certainly true:  if the flow
	downstairs is integrable then the flow
	upstairs is integrable. (See the
paper by Fomenko and Mischenko [FM],
or [A].)  Our assumption is probably
false in general, but we expect
the exceptions to be `pathological'.   

	 We need to say a few words about what
	we mean by  ``the flow downstairs
	being integrable''.  
	 Any Hamiltonian vector field on
	$\g^*$ is necessarily tangent to the orbits
	of the co-adjoint action of $G$ on
	$\g^*$.  The Poisson structure
	induces  a symplectic structure
	on these orbits, sometimes called
	the Kirilov-Kostant-Souriau structure.
	When we say that the flow on $\g^*$
	is `completely integrable' what we mean
	is that  the  Hamiltonian
	 flows restricted to  typical co-adjoint orbits
	are completely integrable in the
	sense just described. By a  typical orbit
	we mean one whose dimension is   maximal, say $2k$. Now
	$2k = n -r $ where $r$ is the rank of the Lie algebra,
	which is the dimension of its maximal Abelian subalgebra.
	In any case, complete integrability
	on the orbit  means that there are k functionally
	independent functions $f_1, \ldots, f_k$
	on the typical orbit which Poisson-commute with each other
	and such that $H$ can be expressed in terms of them.
	We assume that  these functions vary smoothly with
	the orbit. The  typical orbit is defined
	by the vanishing of $r$ functions $C_1, \ldots, C_r$.
	These functions are Casimirs:  they Poisson commute
	with every function on $\g^*$.  Pulled back
	to $T^*G$, and set of $r$ Casimirs forms a  a functional
	basis for the bi-invariant functions on $T^*G$.  In
	the nilpotent case it is known that the Casimirs are
	rational functions (see ....).  
	Thus: to say that the reduced system
	is integrable means that   $H = h(f_1 , \ldots , f_k; C_1,
	\ldots , C_r)$ for some smooth function $h$.
	We say it is {\sc algebraically completely integrable}
	if the $f_i$ and $h$ are rational functions.

	\section{The example}
	\label{sect-ex}
	Take $G$ to be the group of all 4 by 4 lower
	triangular matrices with 1's on
	the diagonal.  Its Lie
	algebra $\n$ is the space of all strictly lower triangular
	matrices:

\[
\pmatrix{
0 & 0 & 0 & 0 \cr
x & 0 & 0 & 0 \cr
u & y & 0 & 0 \cr
w & v & z & 0
}. \label{mat_low}
\]
 	It is generated by the three-dimensional
	subspace $V_1$ consisting of the subdiagonal matrices:
\[
\pmatrix{
0 & 0 & 0 & 0 \cr
x & 0 & 0 & 0 \cr
0 & y & 0 & 0 \cr
0 & 0 & z & 0
}
\]
 	and is coordinatized by $x,y,z$.  
	The functions $x,y,z, u,v ,w$ are linear coordinates
	on the dual of the Lie algebra, and hence
	left-invariant fiber linear
	functions on $T^*G$.  The space
$V_2 = [V_1, V_1]$ is the uv plane,
and $V_3 = [V_1, V_2]$ is the w-axis.
For the inner product on $V_1$
	we take the standard Euclidean one so
that x, y,z correspond to the standard orthogonal
axex.  Thus
	the \sR Hamiltonian is

	\[ H = {1 \over 2} (x^2 + y^2 + z^2). \]

The Poisson structure is the standard Kirillov-Kostant structure
defined by the formula:
\[ \{F,G\}(x)=\langle x,[dF(x),dG(x)]\rangle. \]
It is well known that this structure is nondegenerate 
on the orbits of coadjoint action.
Casimirs, i.e. functions invariant under the coadjoint action, are generated by
$w$ and $uv-yw$, see, for example, \cite{DLNT}.

	\begin{theorem}
	The geodesic flow on $\n ^*$
generated by $H$ is not algebraically
	completely integrable. 
	\end{theorem}

{\bf Proof.}
The Kirillov-Kostant-Souriau Poisson bracket is given by the
 relations:

\[
\matrix{
\{z,u\} = w \cr
\{v,x\} = w \cr
\{y,x\} = u \cr
\{z,y\} = v \cr
}
\]
with all  the other Poisson brackets of the coordinate functions
equalling zero.
 We choose $x$, $z$, $\tu=\dfrac{u}{w}$ and $\tv=\dfrac{v}{w}$ as
 Darboux coordinates
on generic orbit.  ``Generic'' means that   $w=w_0\ne 0$ and $uv-yw=C\ne
0$.  One easily checks that   $x$, $z$, $\tu$ and $\tv$ are
independent coordinates on a generic  orbit and that
$$\matrix{
\{z,\tu\} = 1 \cr
\{\tv,x\} = 1 \cr
\{z,x\}=\{z,\tv\}=\{x,\tu\}=\{\tu,\tv\}=0.\cr
}
$$
 In these coordinates Hamiltonian has a form
$H(x,z,\tu,\tv)=x^2+z^2+\left(\frac{C-uv}{w_0}\right)^2=
x^2+z^2+\left(\frac{C}{w_0}-w_0\tu\tv\right)^2$.

Under the following linear symplectic change of coordinates
\[ \matrix{
x={}^3\sqrt{w_0} \hat x \cr
z={}^3\root \of w_0 \hat z \cr
\tu=\dfrac{\hu}{{}^3\root \of w_0} \cr
\tv=\dfrac{\hv}{{}^3\root\of w_0} \cr
} \]
the Hamiltonian becomes $w_0^\frac{2}{3}({\hat x}^2+{\hat z}^2-2C{}^3\root\of w_0\hu\hv+
\hu^2\hv^2+K)$, where $K$ is some constant depending on $w_0$ and $C$ only.
This Hamiltonian is proportional to a famous Yang-Mills hamiltonian,
which Ziglin proved to be  rationally nonintegrable.  See Ziglin
\cite{Ziglin1},\cite{Ziglin2}.  \qed

\section{Quantization: Extension to the Universal Enveloping Algebra}

The non-integrability in rational functions  of the system considered above 
has some purely algebraic consequences.  

	\begin{lemma}\label{lemma1}
Let $\n$ be the algebra of nilpotent lower triangular $4\times 4$-matrices
and $Pol(\n^*)$ be the algebra of polynomials on its dual space $\n^*$.
Let $\{F,H\}=0$, where $H=\frac{1}{2}(x^2+y^2+z^2)$ and $F\in Pol(\n^*)$.
 Then $F=P(H,w,uv-yw)$. Here $x,y,z,u,v,w$ were defined in the previous
section and $P$ is a polynomial.	
         \end{lemma}

{\bf Proof.}
Taking into account the fact that the dimension of the generic orbit of
the coadjoint representation in $\n^*$ is $4$, we see that if $F$ 
commuted with $H$ but were not  
of the form $P(H,w,uv-yw)$ then the system defined
by $H$ would be completely
integrable. But this contradicts Theorem 1.  \qed

Given any function $f$ in $Pol(\n^*)$
its centralizer with respect to Poisson bracket
always contains the polynomials $F(f, w, uv -yw)$.
For, as mentioned earlier,  
$w$, $uv -yw$ are the {\sc Casimirs} for $\n^*$:
they generate the center of
$Pol(\n^*)$.  Thus the lemma asserts that the
centralizer of H is as small as possible.

We will finish off by proving a similar result for the
universal enveloping algebra $U(\n)$. 
$U(\n)$  can be thought
of as the algebra of left-invariant differential operators
on the Lie group $N$ 
of upper triangular matrices with 1's on the
diagonal (or on certain
homogeneous spaces for N).
It is  
is generated as an algebra
over $\bf R$ by $\n$ (the 1st order
differential operators) and the unit 1 (the identity operator).  Let
$E_{ij}$ be the standard unit matrix with only nonzero $(i,j)$ entry
equal to $1$.  Set $X=E_{21}$, $Y=E_{32}$, $Z=E_{43}$, $U=E_{31}$,
$V=E_{42}$, $W=E_{41}$.  Then  $X, Y, Z, U, V, W$ 
together with 1 generate $U(\n)$.   
 Observe that
this notation is consistent with that used for the elements
$x,y,z,u,v, w$ for $Pol(\n^*)$.  Thus w is a linear function
on $\n^*$, which is to say an element of $\n$.

Let $\HH=\frac{1}{2}({X}^2+{Y}^2+{Z}^2)\in U(\n)$.
It is the ``quantization'' of our
$H \in Pol(\n^*)$.  
By a theorem of Hormander, it is a hypoelliptic
differential operator
which is almost as good as being elliptic.    
 It is well-known that the center of $U(\n)$ is
generated by $W$ and $UV- YW$.  (See
 Dixmier \cite{Dixmier}.)   Consequently if
$R$ is any element of $U(\n)$ then its commutator
algebra {\bf contains} the subalgebra generated
by R, $W$ and $UV- YW$.        

\begin{theorem}\label{th2}
Any element  $F$  in $U(\n)$ which commutes
with $\HH$ is of  the form
$F=P(\HH,W, U V - Y W)$ for some polynomial $P$.
\end{theorem}

There is no ordering problem in defining
the element $P(H, W, UV - YW)$ since $W$ and
$UV-YW$ commute with everything.
 The theorem asserts that
$\HH$ commutes only with those elements which every operator
must commute with.   
It suggests
that   $\HH$, as an operator, should exhibit ``quantum chaos''. 

This theorem is  a special case of a result
which holds for any finite-dimensional Lie
algebra $\g$. The result may be  well-known to experts
but we will present it here in any case.   
 
Let
$U(\g)$ be the  universal enveloping algebra 
of the finite-dimensional Lie algebra $\g$ and 
$Z(G) \subset U(G)$ its center.
  $U(\g)$ is filtered by degree.  
An element is said to have degree less than or equal to k
if it is a sum of monomials 
of the form $X_1 X_2 \ldots X_s$
with the $X_i \in \g$ and $s \le k$.  If
$s = k$ for one of of these monomial terms
then its degree equals k.     
The corresponding graded algebra 
$Gr(U(G))$ is canonically isomorphic
to the algebra $Pol(\g^*)$ of polynomials
on $\g^*$.  The operator bracket
respects the  filtration 
so that it induces a Lie bracket on $Pol(\g^*)$.
This is of course the  
KKS Poisson-bracket $\{, \cdot , \cdot \}$ on $\g^*$.

Let $U(\g)_k$ denote the subspace of elements
of degree k or less.  The quotient map
$\sigma_k: U(\g)_k \to U(\g)_k/U(\g)_{k-1} \cong Pol(\g^*)_k$
takes elements of degree k to homogeneous polynomials
of degree k.   If 
an element $\tilde F \in U(G)$
has  degree  k then we call $\sigma_k(F)$
its principal symbol. If two elements
in $U(\g)$ commute then their principal
symbols must Poisson-commute in $Pol(\g)^*$.
This follows directly from the relation between the
operator and Poisson brackets.

There is a symmetrization map $\phi: Pol(\g^*) \to U(\g)$
which is a kind of 
 inverse to the symbol maps.
It is a linear isomorphism but of course not an
algebra homomorphism.  When restricted to the subspace of 
homogeneous polynomials of
degree k it satisfies $\sigma_k \circ \phi = Id.$

 The center $Z(\g)$ of $U(\g)$
is  finitely generated by elements
$\tilde f_1, \ldots, \tilde f_r$.  These elements
may be chosen so that
their principal symbols 
$f_1, \ldots , f_r$ generate the
 the center of $Pol(\g)$ and so that $\tilde f_i = \phi(f_i)$.
(See Dixmier, or Varadarajan \cite{Var}, Thm. 3.3.8,
p. 183.) Elements of either center are called
{\sc Casimirs}.    (If $\g$ is semi-simple
then the number r of Casimirs is the rank of the Lie algebra.)

\begin{proposition} 
Let   $\tilde H$ be an element of $U(\g)$ of degree
m and   $H = \phi_m(\tilde H)$ its principal symbol.   
Suppose that the commutator algebra of $H$ in
$Pol(\g^*)$ is generated by the Casimirs
$f_1, \ldots , f_r$
together with $H$.  Then the commutator algebra of
$\tilde H$ in $U(\g)$ is generated
by the Casimirs $\tilde f_1, \ldots , \tilde f_m$ 
together with $\tilde H$.
\end{proposition}

In other words, if the commutator of the
principal symbol is as small as possible, the
same is true for its quantization $\tilde H$.  

{\bf Proof.}  Suppose $\tilde F$ commute with
$\tilde H$.  Let k denote the
degree of $\tilde F$ and set  $F = \sigma_k(\tilde F)$.
As discussed above, F 
Poisson-commutes with H. 
By hypothesis $F = p(f_1, \ldots, f_r, H)$
for some polynomial $p$.  
Since the $\tilde f_i$ are in the center of
$U(G)$, the element
$\tilde p = p(\tilde f_1, \ldots \tilde f_r, \tilde H)$
is a well-defined element
of $U(\g)$, independent of the  ordering of
the factors $\tilde f_i, \tilde H$.
Clearly $\tilde p$ commutes with
$\tilde H$.  
Moreover $\sigma_k(\tilde p) = F$ so that
$\sigma_k(\tilde F - \tilde p) = 0$.
 It follows that $\tilde F - \tilde p = \tilde F_2$
is an element whose degree is $k-1$ or less
which commutes with $\tilde H$.

Repeating this argument with $\tilde F_2$ in
place of $\tilde F$ we
obtain a polynomial $p_2$
in $f_1 , \ldots, f_r,H$ and a corresponding
element $\tilde p_2$ in $U(\g)$ such that
$\tilde F - (\tilde p + \tilde p_2)$ commutes
with $\tilde H$ and has degree $k-2$ or less.
Continuing in this fashion we eventually descend
to degree 0, in which case :
$$\tilde F = \tilde p + \tilde p_2 + \ldots
+  \tilde p_{k-1}$$
is a polynomial in the $\tilde f_i, \tilde H$
as claimed.
QED

\end{document}